\documentclass[12pt]{article}
       \usepackage{amsmath,amsthm}
       \usepackage{amssymb,amsbsy}
       \usepackage{amsfonts}

\usepackage{times}
\usepackage{cite}

\begin{document}

\author{
{\large {\bf  \framebox{W. Marcinek}
}}\\{\it Institute of Theoretical Physics, University of Wroc\l aw,}\\Pl. Maxa Borna 9,
50-204 Wroc{\l}aw, Poland
}

\title{{\bf QUANTUM YANG-BAXTER
EQUATION, ARBITRARY GRADINGS\\ AND EXCHANGE HECKE BRAIDINGS}
\thanks{This is the last finished paper by W\l adys\l aw Marcinek
who unexpectedly passed away on June 9, 2003. Prepared for publication by Steven Duplij.}
\thanks{Published in the W. Marcinek's Memorial volume of the
Kharkov University Journal  Vestnik KNU, ser. "Nuclei, particles and fields".
- 2003. - V. 601. - N 2(22). - p. 60--66.}
}
\thispagestyle{empty}
 \date{May 15, 2003}
\maketitle
\begin{abstract}
Generalization of the quantum Yang-Baxter equation solutions to an
arbitrary grading is studied. The noncommutative differential calculi
corresponding to such solutions is considered. The connection with the
ordinary and supersymmetric solutions of the quantum Yang-Baxter equation is
given.
\smallskip


\end{abstract}
\newpage
It is well known that the supersymmetry and the superspace formalism has been
developed in the last 20 - 30 years in the context of Fermi-Bose
correspondence, supergravity and string theory. A few generalizations of
supersymmetry are known. The so-called color generalization of supersymmetry
transformation has been introduced into physics by Rittenberg and Lukierski
\cite{4}, Rittenberg and Wyler \cite{rw,rwu}. The corresponding algebraic
structures have been developed by several authors: Scheunert \cite{5},
Kobayashi and Nagamachi \cite{kn}, Trostel \cite{tro}, Kwa\'{s}niewski
\cite{kk}, Marcinek \cite{6,WM2,WM3,10} and many others. Similar topics have
been also considered by Matthes \cite{RM}. Next, the formal $S$-group of
transformation as a generalization of supersymmetry, where $S$ is a triangular
Yang-Baxter operator (sometime called a symmetry), has been introduced by
Gurevich \cite{Gur,G,3}. The mathematical structures corresponding to such
transformations has been studied by Lyubashenko \cite{9}, and further
developed by others, see \cite{man,10,WM5} for example. The generalization for
arbitrary quasitriangular solution (braidings) of the quantum Yang-Baxter
equation has been given by Majid \cite{8,Maj,bm,qm}. In all of these above
formalisms the tensor product of spaces is modified.

On the other hand the $q$-deformed structures and related quantum groups based
on the usual (bosonic) tensor product, have been lately under intensive study
in both, mathematics and theoretical physics. The corresponding formalism has
been developed by several authors, Drinfeld, Jimbo, Woronowicz and many
others. The well-known in theoretical physics quantum groups correspond to the
$q$-deformation of universal enveloping algebras of simple Lie algebras. It is
interesting that they admits the so-called universal $R$-matrix which also
provide a solution of quantum Yang-Baxter equation. For example the $R$-matrix
corresponding to quantum deformations of the groups of the series $A_{n-1}$
satisfy the so-called Hecke condition. Note that we have in fact two attempts
to quantum deformations. In the first approach only the tensor product is
deformed, the algebra structure remains undeformed. In physical interpretation
this case correspond to nonstandard statistics of particles. In the second
approach the situation is opposite, the algebra structure is deformed, but the
tensor product - not. In these two approach the concept of $R$-matrix play an
important role. It is interesting that these two approach can be unified, see
Majid \cite{8} and Hlavaty \cite{HL}. In this way we obtain formalism with
deformed tensor product and algebra structure. The corresponding formalism can
be understood as a quantum deformations with braid group statistics. The
quantum superspace studied by Ilinski and Uzdin \cite{iuz} can be considered
as an example of such structures. Graded structures studied previously by the
author in \cite{WM3} have been given as an another example for such structures.

In this paper we study some new graded structures which unify the color
symmetry and quantum deformations. We introduce the concept of exchange
braidings and color $R$-matrix. We study in details the exchange braiding and
color $R$-matrix corresponding for quantum deformations with Hecke condition.
A few examples for such $R$-matrices are given. They can be understood as a
color generalizations of $q$-deformations corresponding to the Lie (super-)
groups of the series $A_{n-1}$, \cite{12,13}. The relation of our generalized
graded structures with the standard (ungraded) or supersymmetric formalism is
investigated. As a result we obtain that the odd (resp. even) color Hecke
$R$-matrix can be reduced to a certain tensor product of $Z_{2}$-graded (resp.
ungraded) solution of the quantum Yang-Baxter equation with Hecke condition
and even color $R$-matrix. The generalization of Wess-Zumino \cite{17}
noncommutative differential calculi on a quantum plane is considered.

\section{EXCHANGE BRAIDINGS}

Let $E$ be an arbitrary vector space over a field $k$. A linear operator
$R:E\otimes E\longrightarrow E\otimes E$ such that%
\begin{equation}%
\begin{array}
[c]{c}%
R_{12}\ R_{23}\ R_{12}=\ R_{23}\ R_{12}\ R_{23}\label{ybe}%
\end{array}
\end{equation}
is said to be \textit{a Yang-Baxter operator or a quasisymmetry or a braiding}
on $E$, \cite{man,2,3}, where $R_{12}=R\otimes id$, and $R_{23}=id\otimes R$.
If in addition
\begin{equation}%
\begin{array}
[c]{c}%
R^{2}=id+(q-q^{-1})R,
\end{array}
\end{equation}
then $R$ is said to be a Hecke symmetry on $E$, or a Hecke braiding \cite{3}.

Now let $E$ be a $\Gamma$-graded vector space over a field of complex numbers
$k\equiv\mathbf{C}$, i.e. we have a direct sum decomposition $E=\bigoplus
_{\alpha\in\Gamma}E_{\alpha}$. If the grading group $\Gamma$ has a nontrivial
subgroup $\Gamma_{0}$ of index $2$ (i. e. the quotient group $\Gamma
/\Gamma_{0}$ is isomorphic to the group $Z_{2}$), then the above gradation is
said to be \textit{odd}, in the opposite case it is called \textit{even}. We
denote by $\pi:\Gamma\longrightarrow\Gamma/\Gamma_{0}$ the corresponding
quotient map. Assume that we have a homogeneous basis $\{e_{i}:i\in I\}$ of
the space $E$, $I=\{1,2,\ldots,N\}$, $N=\mathrm{\dim}E$. This means that
$e_{i}$ is the $i$-th basis vector of $E$ of grade $p(i)$, $p:I\longrightarrow
\Gamma$ is a gradation mapping, i.e. $p(i)=\alpha$ if $e_{i}\in E_{\alpha}$.
We introduce here a double index notation. We denote by $e_{i_{\alpha},\alpha
}$ the $i_{\alpha}$-th basis vector of $E$. This means that for fixed $\alpha$
the set $\{e_{i_{\alpha},\alpha}:i\in I_{\alpha}\}$ form a basis for
$E_{\alpha}$, where $I_{\alpha}:=\{i\in I:p(i)=\alpha\}$. Observe that
$\pi(\alpha)=+1$ or $0$. The value of $\varepsilon_{\alpha}:=(-1)^{\pi
(\alpha)}$ is said to be \textit{a parity} of $\alpha$. Note that if
$\varepsilon_{\alpha}=+1$ for every $\alpha\in G$, then the gradation is even.
If there is an $\alpha\in\Gamma$ such that $\varepsilon_{\alpha}=-1$, then it
is odd.

Let $\varepsilon$ be a commutation factor on $\Gamma$, i.e. a mapping
$\varepsilon:\Gamma\times\Gamma\longrightarrow C\diagdown\{0\}$ such that
$\varepsilon(\alpha,\beta)\varepsilon(\beta,\alpha)=1$ and $\varepsilon
(\alpha+\beta,\gamma)=\varepsilon(\alpha,\gamma)\varepsilon(\beta,\gamma)$ for
every $\alpha,\beta,\gamma\in\Gamma$ \cite{5}. Note that $\varepsilon
(\alpha,\alpha)=+1$, or $-1$ for every $\alpha\in\Gamma$. If $\varepsilon
(\alpha,\alpha)=+1$ for every $\alpha\in\Gamma$, then $\varepsilon$ is said to
be even commutation factor. If there is an $\alpha\in\Gamma$ such that
$\varepsilon(\alpha,\alpha)=-1$, then $\varepsilon$ is said to be odd
\cite{kn}. The set $\Gamma_{0}=\{\alpha\in\Gamma:\varepsilon(\alpha
,\alpha)=1\}$ is obviously a subgroup of $\Gamma$ of index at most $2$. If the
factor $\varepsilon$ is odd, then the subgroup $\Gamma_{0}$ is nontrivial and
the corresponding gradation is odd. Observe that here $\varepsilon_{\alpha
}\equiv\varepsilon(\alpha,\alpha)$. For a given $\Gamma$-graded vector space
$E$ equipped with a homogeneous basis $\{e_{i_{\alpha},\alpha}:i\in I_{\alpha
},\alpha\in\Gamma\}$ we introduce the concept of exchange operators.

A linear and homogeneous mapping $R:E\otimes E\longrightarrow E\otimes E$,
where $E$ is a $\Gamma$-graded vector space, which satisfies the Yang-Baxter
equation (\ref{ybe}) on $E$ is said to be \textit{a graded Yang-Baxter
operator }on $E$. It is obvious that for every graded Yang-Baxter operator $R$
we have the following decomposition
\begin{equation}%
\begin{array}
[c]{l}%
R:=\bigoplus\limits_{\alpha,\beta\in\Gamma}\ R(\alpha,\beta),\label{pos}%
\end{array}
\end{equation}
where $R(\alpha,\beta):E_{\alpha}\otimes E_{\beta}\longrightarrow E_{\beta
}\otimes E_{\alpha}$ is the homogeneous component of $R$ with respect to the
$\Gamma\times\Gamma$-gradation of the tensor product $E\otimes E$. Every such
component $R(\alpha,\beta)$ of $R$ is said to be \textit{a color exchange
operator} or (\textit{a color R-matrix}) on $E$. The following result is
obvious: the quantum Yang-Baxter equation (\ref{ybe}) for a graded Yang-Baxter
operator $R$ on $E$ is equivalent to the following family of equations
\begin{equation}%
\begin{array}
[c]{c}%
R_{12}(\beta,\gamma)R_{23}(\alpha,\gamma)R_{12}(\alpha,\beta)=R_{23}%
(\alpha,\beta)R_{12}(\alpha,\gamma)R_{23}(\beta,\gamma),\label{cyb}%
\end{array}
\end{equation}
defined for the family of exchange operators $\{R(\alpha,\beta):\alpha
,\beta\in\Gamma\}$, on $E_{\alpha}\otimes E_{\beta}\otimes E_{\gamma}$ for
every $\alpha,\beta,\gamma\in\Gamma$. The family of relations (\ref{cyb}) is
called \textit{a color exchange quantum Yang-Baxter equation}. The above
result means that every graded Yang-Baxter operator $R$ on $E$ can be
determined by a family of color exchange operators $\{R(\alpha,\beta
):\alpha,\beta\in\Gamma\}$ defined for a family of spaces $\{E_{\alpha}%
:\alpha\in\Gamma\}$, where $E_{\alpha}$ is the $\alpha$-th component of $E$
with respect to a given $\Gamma$-gradation. The family $\{R(\alpha,\beta)\}$
is sometime said to be \textit{a color exchange} braidings (or symmetries) for
$\{E_{\alpha}:\alpha\in\Gamma\}$. Let us consider an example: a family of
operators $\{(R_{\varepsilon})(\alpha,\beta):\alpha,\beta\in\Gamma\}$ given by
the following relation
\begin{equation}%
\begin{array}
[c]{c}%
R_{\varepsilon}(\alpha,\beta)(e_{i_{\alpha},\alpha}\otimes e_{j_{\beta},\beta
}):=\varepsilon(\alpha,\beta)\ e_{j_{\beta},\beta}\otimes e_{i_{\alpha}%
,\alpha},\label{col}%
\end{array}
\end{equation}
where $\varepsilon$ is a commutation factor on the grading group $\Gamma$, is
a family of exchange operators for $\{E_{\alpha}:\alpha\in\Gamma\}$. This
family define a graded Yang-Baxter operator $R_{\varepsilon}$ on $E$. The
operator $R_{\varepsilon}$ is said to be color.

If the commutation factor $\varepsilon$ on $\Gamma$ is even (resp. odd), then
these above color operators or symmetries are said to be even or even (resp.
odd). The above graded Yang-Baxter operator is diagonal with respect to both
classes of indices: external and internal. Now let us consider the case in
which a given graded Yang-Baxter operator is diagonal only in external
(grading) indices. To this goal we introduce some special grading of
$E_{\Gamma}$. Let $E_{\Gamma}$ be a $\Gamma$-graded space such that every
component $E_{\alpha}$ of $E$ with respect to a given $\Gamma$-gradation has
the same dimension equal to $k$, i.e. $\dim(E_{\alpha})=k$ for every
$\alpha\in\Gamma$. This means that all these components $E_{\alpha}$ of
$E_{\Gamma}$ are isomorphic one to other. In other words, there exist a space
$E$ such that each $E_{\alpha}$ is isomorphic to $E$. In this way every
$E_{\alpha}$ can be understand as a copy of the space $E$ and the $\Gamma
$-graded vector space $E_{\Gamma}$ which is direct sum of all spaces
$E_{\alpha}$ can be identified with the product $E\times\Gamma$. Hence we have
here an analogy with the topological concept of covering spaces in a certain
algebraic sense. The space $E_{\Gamma}$ play in fact a role of algebraic
covering of the "base" space $E$ with a discrete "fiber" $\Gamma$. Let us
denote by $e_{i}$ the $i$-th basis vector of $E$, $i=1,\ldots,k$. The one-to
one correspondence between spaces $E_{\Gamma}$ and $E\times\Gamma$ is given by
$e_{i_{\alpha},\alpha}\leftrightarrow(e_{i},\alpha)$, where $i_{\alpha
}=1,\ldots,N$, $i=1,\ldots,k$, $\alpha\in\Gamma$. The tensor product is
$E_{\Gamma}\otimes E_{\Gamma}$.

Observe that for every element $e_{i_{\alpha},\alpha}$ of $E_{\alpha}$ there
is a corresponding element $e_{i}$ of $E$. We define a set of components
$R(\alpha,\beta)$ for an operator $R$ on a $\Gamma$-graded space $E$ by the
following formula
\begin{gather}
R(\alpha,\beta)=q^{1-2\pi(\alpha)}\ \varepsilon(\alpha,\alpha)\ \Sigma_{i\in
I_{\alpha}}\ e_{i,\alpha\ i,\alpha}\otimes e_{i,\alpha\ i,\alpha}\nonumber\\
+\Sigma_{i\in I_{\alpha},j\in I_{\beta},i\neq j}\ \varepsilon(\alpha
,\beta)\ e_{i,\alpha\ j,\beta}\otimes e_{j,\beta\ i,\alpha}+\varepsilon
(\alpha,\alpha)(q-q^{-1})\nonumber \\ \times\Sigma_{i\in I_{\alpha},j\in I_{\beta}%
,i<j}\ e_{i,\alpha\ i,\alpha}\otimes e_{j,\beta\ j,\beta}, \label{per}%
\end{gather}
where $\varepsilon$ is a commutation factor on $\Gamma$, $q$ is a deformation
parameter, $q\in C\diagdown\{0\}$, and $(e_{i,\alpha\ j,\beta})_{k,\gamma
\ l,\delta}=\delta_{i,\alpha\ k,\delta}\delta_{j,\beta\ l,\alpha}$.

Note that the corresponding family of mappings $R(\alpha,\beta):E_{\alpha
}\otimes E_{\beta}\longrightarrow E_{\beta}\otimes E_{\alpha}$ is given
explicit by the equation
\begin{equation}
R(\alpha,\beta)(e_{i,\alpha}\otimes e_{j,\beta})=\left\{
\begin{array}
[c]{lll}%
q^{1-2\pi(\alpha)}\varepsilon(\alpha,\alpha)e_{i,\alpha}\otimes e_{i,\alpha
}, & \text{for} & i=j,\\
\varepsilon(\alpha,\alpha)(q-q^{-1})e_{i,\beta}\otimes e_{j,\alpha
}+\varepsilon(\alpha,\beta)e_{j,\beta}\otimes e_{i,\alpha}, & \text{for} &
i<j,\\
\varepsilon(\alpha,\beta)e_{j,\beta}\otimes e_{i,\alpha}, & \text{for} & i>j.
\end{array}
\right.  \label{sys}%
\end{equation}
If the commutation factor $\varepsilon$ on $\Gamma$ is even, then the family
(\ref{ope}) is said to be \textit{even color Hecke braidings}. If
$\varepsilon$ is odd, then this family is also said to be \textit{odd}. Next
we use the notation: $\varepsilon_{\alpha}:=\varepsilon(\alpha,\beta
)\equiv(-1)^{\pi(\alpha)}$.

For the family of operators $\{R(\alpha,\beta)\}$ given by the equation
(\ref{ope}) we have \textbf{(i)} the Yang-Baxter equation (\ref{cyb}), and
\textbf{(ii)} the corresponding Hecke condition
\begin{equation}%
\begin{array}
[c]{c}%
R^{2}(\alpha,\beta)=id+(-1)^{\pi(\alpha)}(q-q^{-1})R(\alpha,\beta).
\end{array}
\end{equation}

To show this first we prove the Yang-Baxter equation (\ref{cyb}) for the
operation (\ref{ope}). For $i<j<k$ we calculate
\begin{align}
&R_{12}(\alpha,\beta)(e_{i,\alpha}\otimes e_{j,\beta}\otimes x_{k,\gamma
})=R(\alpha,\beta)(e_{i,\alpha}\otimes e_{j,\beta})\otimes e_{j,\gamma
}\nonumber\\&=\varepsilon_{i,\alpha}(q-q^{-1})e_{i,\alpha}\otimes e_{j,\beta}\otimes
e_{j,\gamma}+\varepsilon(\alpha,\beta)e_{j,\beta}\otimes e_{i,\alpha}\otimes
e_{j,\gamma},
\end{align}
and
\[%
\begin{array}
[c]{l}%
R_{23}(\alpha,\gamma)R_{12}(\alpha,\beta)(e_{i,\alpha}\otimes e_{j,\beta
}\otimes e_{j,\gamma})=\varepsilon_{i,\alpha}(q-q^{-1})e_{i,\alpha}\otimes
e_{j,\beta}\otimes e_{j,\gamma}\nonumber\\+\varepsilon_{i,\alpha}(q-q^{-1})\varepsilon
(\beta,\gamma)e_{i,\alpha}\otimes e_{j,\gamma}\otimes e_{j,\beta}\nonumber
+\varepsilon_{i,\alpha}(q-q^{-1})\varepsilon(\alpha,\beta)e_{j,\beta}\otimes
e_{i,\alpha}\otimes e_{j,\gamma}\\+\varepsilon(\alpha,\beta+\gamma)e_{j,\beta
}\otimes e_{j,\gamma}\otimes e_{i,\alpha}.
\end{array}
\]
Hence for the left hand side of the Yang-Baxter equation (\ref{cyb}) we
obtain
\[%
\begin{array}
[c]{l}%
R_{12}R_{23}R_{12}(e_{i,\alpha}\otimes e_{j,\beta}\otimes e_{j,\gamma
})=[\varepsilon_{i,\alpha}(q-q^{-1})^{3}+\varepsilon_{i,\alpha}(q-q^{-1}%
)]e_{i,\alpha}\otimes e_{j,\beta}\otimes e_{j,\gamma}\\+(q-q^{-1})^{2}%
\varepsilon(\alpha,\beta)e_{j,\beta}\otimes e_{i,\alpha}\otimes e_{j,\gamma}
+(q-q^{-1})^{2}\varepsilon(\beta,\gamma)e_{i,\alpha}\otimes e_{j,\gamma
}\otimes e_{j,\beta}\\+\varepsilon_{i,\alpha}(q-q^{-1})\varepsilon(\alpha
+\beta,\gamma)e_{j,\gamma}\otimes e_{i,\alpha}\otimes e_{j,\beta}%
+\varepsilon_{i,\alpha}(q-q^{-1})\varepsilon(\alpha,\beta+\gamma)e_{j,\beta
}\otimes e_{j,\gamma}\otimes e_{i,\alpha}\\
+\varepsilon(\alpha,\beta)\varepsilon(\alpha,\gamma)\varepsilon(\beta
,\gamma)e_{j,\gamma}\otimes e_{j,\beta}\otimes e_{i,\alpha},
\end{array}
\]
and for the right-hand side
\[%
\begin{array}
[c]{l}%
R_{23}R_{12}R_{23}(e_{i,\alpha}\otimes e_{j,\beta}\otimes e_{j,\gamma
})=[\varepsilon_{i,\alpha}(q-q^{-1})^{3}+\varepsilon_{i,\alpha}(q-q^{-1}%
)]e_{i,\alpha}\otimes e_{j,\beta}\otimes e_{j,\gamma}\\+(q-q^{-1})^{2}%
\varepsilon(\beta,\gamma)e_{i,\alpha}\otimes e_{j,\gamma}\otimes e_{j,\beta}
+(q-q^{-1})^{2}\varepsilon(\alpha,\beta)e_{j,\beta}\otimes e_{i,\alpha}\otimes
e_{j,\gamma}\\+\varepsilon_{i,\alpha}(q-q^{-1})\varepsilon(\alpha,\beta
+\gamma)e_{j,\beta}\otimes e_{j,\gamma}\otimes e_{i,\alpha}+\varepsilon
_{i,\alpha}(q-q^{-1})\varepsilon(\alpha+\beta,\gamma)e_{j,\gamma}\otimes
e_{i,\alpha}\otimes e_{j,\beta}\\
+\varepsilon(\alpha,\beta)\varepsilon(\beta,\gamma)\varepsilon(\alpha
,\gamma)e_{j,\gamma}\otimes e_{j,\beta}\otimes e_{i,\alpha}.
\end{array}
\]
The proof in the remaining cases is obvious. We have for $i<j$
\[%
\begin{array}
[c]{l}%
R^{2}(\alpha,\beta)(e_{i,\alpha}\otimes e_{j,\beta})=R(\beta,\alpha
)((q-q^{-1})e_{i,\alpha}\otimes e_{j,\beta}+\varepsilon(\alpha,\beta
)e_{j,\beta}\otimes e_{i,\alpha})\\
=\varepsilon_{i,\alpha}(q-q^{-1})^{2}x_{i,\alpha}\otimes x_{j,\beta
}+\varepsilon_{i,\alpha}(q-q^{-1})\varepsilon(\alpha,\beta)e_{j,\beta}\otimes
e_{i,\alpha}+e_{i,\alpha}\otimes e_{j,\beta}\\
=(id+(q-q^{-1})R)(e_{i,\alpha}\otimes e_{j,\beta}).
\end{array}
\]

The operation defined by the equation (\ref{ope}) is called \textit{a color
Hecke R-matrix} and the corresponding family of exchange mappings is said to
be \textit{a color Hecke braiding} on a $\Gamma$-graded space $E$
corresponding to the commutation factor $\varepsilon$.

Let us consider some particular examples. We define a operator $R$ on a
$\Gamma$-graded space $E$ by the following formula
\begin{align}
&R    =\bigoplus\limits_{\alpha,\beta\in\Gamma}(q^{1-2\pi(\alpha
)}\ \varepsilon(\alpha,\alpha)\ \Sigma_{i\in I_{\alpha}}\ e_{i,\alpha
\ i,\alpha}\otimes e_{i,\alpha\ i,\alpha}\nonumber\\
&  +\Sigma_{i\in I_{\alpha},j\in I_{\beta},i\neq j}\ \varepsilon(\alpha
,\beta)\ e_{j,\beta\ i,\alpha}\otimes e_{i,\alpha\ j,\beta}+\varepsilon
(\alpha,\alpha)(q-q^{-1})\nonumber\\&\times \Sigma_{i\in I_{\alpha},j\in I_{\beta}%
,i<j}\ e_{i,\alpha\ i,\alpha}\otimes e_{j,\beta\ j,\beta}), \label{ope}%
\end{align}
where $\varepsilon$ is a commutation factor on $\Gamma$, $q$ is a deformation
parameter, $q\in C\diagdown\{0\}$, and $(e_{i,\alpha\ j,\beta})_{k,\gamma
\ l,\delta}=\delta_{i,\alpha\ k,\delta}\delta_{j,\beta\ l,\alpha}$. Note that
the corresponding family of mappings $R(\alpha,\beta):E_{\alpha}\otimes
E_{\beta}\longrightarrow E_{\beta}\otimes E_{\alpha}$ is given explicit by the
equation
\begin{equation}
R(\alpha,\beta)(e_{i,\alpha}\otimes e_{j,\beta})=\left\{
\begin{array}
[c]{lll}%
q^{1-2\pi(\alpha)}\varepsilon(\alpha,\alpha)e_{i,\alpha}\otimes e_{i,\alpha
}, & \text{for} & i=j,\\
\varepsilon(\alpha,\alpha)(q-q^{-1})e_{i,\beta}\otimes e_{j,\alpha
}+\varepsilon(\alpha,\beta)e_{j,\beta}\otimes e_{i,\alpha}, & \text{for} &
i<j,\\
\varepsilon(\alpha,\beta)e_{j,\beta}\otimes e_{i,\alpha}, & \text{for} & i>j.
\end{array}
\right.  \label{esy}%
\end{equation}

Let us consider a $\Gamma$-graded vector space $E$ equipped with a homogeneous
basis ${e_{i}:i\in I}$ and a grading map $p:I\longrightarrow\Gamma$. Let us
assume for simplicity that $\dim(E_{\alpha})=2$ for every every $\alpha
\in\Gamma$ and $\varepsilon_{i}=1$ for $i=1,2$. We use the notation
$e_{1,\alpha}:=x_{\alpha}$, $e_{2,\beta}:=y_{\beta}$.

In this case the color
Hecke $R$-matrix (\ref{ope}) becomes
\begin{align}
&R(\alpha,\beta)=\\&\left[
\begin{array}
[c]{llll}%
q^{1-2\pi(\alpha)}\varepsilon(\alpha,\beta)\delta_{\alpha}^{\delta}%
\delta_{\beta}^{\gamma}\delta_{\alpha\beta} & 0 & 0 & 0\\
0 & (q-q^{-1})\delta_{\alpha}^{\gamma}\delta_{\beta}^{\delta} & \varepsilon
(\alpha,\beta)\delta_{\alpha}^{\delta}\delta_{\beta}^{\gamma} & 0\\
0 & \varepsilon(\beta,\alpha)\delta_{\alpha}^{\delta}\delta_{\beta}^{\gamma} &
0 & 0\\
0 & 0 & 0 & q^{1-2\pi(\alpha)}\varepsilon(\alpha,\beta)\delta_{\alpha}%
^{\delta}\delta_{\beta}^{\gamma}\delta_{\alpha\beta}%
\end{array}
\right]  \nonumber,
\end{align}
where $\alpha,\beta,\gamma,\delta\in\Gamma$.
The corresponding exchange
braiding has the following form
\begin{equation}%
\begin{array}
[c]{l}%
R(x_{\alpha}\otimes x_{\alpha})=q^{1-2\pi(\alpha)}\varepsilon(\alpha
,\alpha)x_{\alpha}\otimes x_{\alpha},\\
R(x_{\alpha}\otimes y_{\beta})=(q-q^{-1})x_{\alpha}\otimes y_{\beta
}+\varepsilon(\alpha,\beta)y_{\beta}\otimes x_{\alpha},\\
R(y_{\beta}\otimes x_{\alpha})=\varepsilon(\beta,\alpha)x_{\alpha}\otimes
y_{\beta},\\
R(y_{\beta}\otimes y_{\beta})=q^{1-2\pi(\beta)}\varepsilon(\beta
,\beta)y_{\beta}\otimes y_{\beta}.\label{25}%
\end{array}
\end{equation}

Let $I=\{1,2,3,4\}$, $\Gamma=Z_{2}$, and the grading map is given by:
$p(1)=p(2)=0,p(3)=p(4)=1$. We have the following supercommutation factor on
$Z_{2}$
\begin{equation}%
\begin{array}
[c]{c}%
\varepsilon(\alpha,\beta)=(-1)^{\alpha\beta}\label{26}%
\end{array}
\end{equation}
for $\alpha,\beta\in Z_{2}$. Here the color $R$-matrix (\ref{ope}) takes the
following form
\[
R=\left[
\begin{array}
[c]{cccc}%
q & 0 & 0 & 0\\
0 & (q-q^{-1}) & 1 & 0\\
0 & 1 & 0 & 0\\
0 & 0 & 0 & -q^{-1}%
\end{array}
\right]  .
\]
This is the well-known (braid) $R$-matrix corresponding to quantum group
$gl(1|1)$, see e.g. \cite{11}.

Let $E$ be a $\Gamma$-graded vector space, where $\Gamma$ be an Abelian group
generated by $\xi_{i}$, $i=1,...,k$, and we have $\mathrm{\dim}E_{\alpha}=1$
for every $\alpha\in\Gamma$. We assume here that the grading map
$p:I\longrightarrow\Gamma$ is given by
\begin{equation}%
\begin{array}
[c]{c}%
p(i)=\xi_{i}.\label{27}%
\end{array}
\end{equation}
This is the so-called standard gradation \cite{10}. We use the notation
$\varepsilon(\xi_{i},\xi_{j})=\varepsilon_{ij}$. Let $\Gamma=Z\oplus...\oplus
Z$ ($k$-summands). Here we have \cite{7,10}
\begin{equation}%
\begin{array}
[c]{c}%
\varepsilon_{ij}=(-1)^{s_{ij}}z^{a_{ij}},\label{28}%
\end{array}
\end{equation}
where $s_{ij}=s_{ji}$ and $a_{ij}=-a_{ji}$ are integer valued matrices, $z\in
C-\{0\}$ is an another deformation parameter, and $\xi_{i}=(0...1...0)$, ($1$
on the $i$-th place). Here the subgroup $\Gamma_{0}$ is generated by the set
$\{\xi_{i}:s_{ii}=0\}$. In this case the color Hecke $R$-matrix (\ref{ope})
can be given in the following form
\begin{equation}
R(\xi_{i},\xi_{j})=q^{1-2s_{ii}}(-1)^{s_{ii}}\Sigma_{i\in I}e_{ii}\otimes
e_{ii}+\Sigma_{i,j\in I,i\neq j}\varepsilon_{ij}e_{ij}\otimes e_{ji}%
+(q+q^{-1})\Sigma_{i,j\in I,i<j}e_{ii}\otimes e_{jj}, \label{29}%
\end{equation}
here $\varepsilon_{ij}$ is given by the equation (\ref{28}). Note that in the
case of $s_{ii}=0$ for all $i\in I$, the the above $R$-matrix is similar the
so-called multiparameter $R$-matrix, see \cite{15,16}. The multiparametric
$R$-matrices are in fact the color Hecke $R$-matrices with $Z^{k}$ -gradings
defined by (\ref{27}) and (\ref{28}), where $s_{ii}=0$ for every $i\in I$.

If $\Gamma=Z_{N}\oplus...\oplus Z_{N}$ ($k$-summands), $N>2$, and
\begin{equation}%
\begin{array}
[c]{c}%
\varepsilon_{ij}=exp((2\pi i/N)\Omega_{ij}),\label{30}%
\end{array}
\end{equation}
where $\Omega_{ij}=-\Omega_{ji}$ is an integer-valued matrix, then the color
Hecke $R$-matrix (\ref{ope}) takes the following form
\begin{equation}
R(\xi_{i},\xi_{j})=q\Sigma_{i\in}e_{ii}\otimes e_{ii}+\Sigma_{i,j\in I,i\neq
j}exp((2\pi i/N)\Omega_{ij})e_{ji}\otimes e_{ij}+(q+q^{-1})\Sigma_{i,j\in
I,i<j}e_{ii}\otimes e_{jj}. \label{31}%
\end{equation}

The above $R$-matrix can be understood as the $R$-matrix corresponding to the
deformation of the anyonic statistics \cite{18}.

\section{REDUCTION OF COLOR-HECKE BRAIDING}

Let $E$ be a $\Gamma$-graded vector space equipped with a homogeneous basis
$\{e_{i,\alpha}:i\in I_{\alpha};\alpha\in\Gamma\}$. We assume here that an odd
color Hecke $R$-matrix on $E$ corresponding to a commutation factor
$\varepsilon$ on $\Gamma$ is given. It follows immediately from our
definitions that the subgroup $\Gamma_{0}$ of $\Gamma$ is nontrivial, i.e. the
quotient $\Gamma/\Gamma_{0}$ is isomorphic to the $Z_{2}$ group. Let us denote
by $\tilde{E}$ the same vector space $E$ but equipped with $Z_{2}$ gradation
defined by the following formulae
\begin{equation}%
\begin{array}
[c]{l}%
\tilde{E}=\tilde{E}_{0}\oplus\tilde{E}_{1},\;\;\tilde{E}_{0}=\bigoplus
_{\alpha\in\Gamma_{0}}E_{\alpha},\;\;\tilde{E}_{1}=\bigoplus_{\alpha\in
\Gamma_{1}}E_{\alpha},\label{zed}%
\end{array}
\end{equation}
where $\Gamma_{1}:=\{\alpha\in\Gamma:\varepsilon(\alpha,\alpha)=-1\}$. The
corresponding $Z_{2}$-grading mapping $\delta:I\longrightarrow Z_{2}$ is given
here by the following formula
\begin{equation}%
\begin{array}
[c]{c}%
\delta=\pi\circ p.\label{det}%
\end{array}
\end{equation}

We denote by $\tilde{e}_{i,\pi(\alpha)}$ the same element $e_{i,\alpha}$ of
$E$ but equipped with the above $Z_{2}$-gradation instead of the initial
$\Gamma$-gradation. We introduce a $\Gamma$-graded vector space $E^{\prime}$
such that
\begin{equation}%
\begin{array}
[c]{c}%
E^{\prime}=\bigoplus_{\alpha\in\Gamma}E_{\alpha}^{\prime},\;\;\dim E_{\alpha
}^{\prime}=1\hspace{0.3cm}for\hspace{0.3cm}every\hspace{0.3cm}\alpha\in
\Gamma,\label{dwa}%
\end{array}
\end{equation}
and $E^{\prime}$ is spanned by element $\rho_{\alpha}$, $\alpha\in\Gamma$. We
denote by $\bar{E}$ the subspace of the tensor product $\tilde{E}\otimes
E^{\prime}$ spanned by elements of the form $\tilde{e}_{i,\alpha}\otimes
\rho_{\alpha}$. Observe that the space $\bar{E}$ is a $\Gamma$-graded vector
space isomorphic to $E$. The space $\bar{E}$ is also $Z_{2}$-graded. The
$Z_{2}$-gradation is given by the grading map (\ref{det}). It is known that
the commutation factor $\varepsilon$ on $\Gamma$ can be given in the following
form
\begin{equation}%
\begin{array}
[c]{c}%
\varepsilon(\alpha,\beta)=(-1)^{\pi(\alpha)\pi(\beta)}\delta(\alpha
,\beta),\label{rco}%
\end{array}
\end{equation}
where $\alpha,\beta\in\Gamma$, and $\delta(\alpha,\beta):=\sigma(\alpha
,\beta)/\sigma(\beta,\alpha)$, and $\sigma$ is a certain $2$-cocycle on
$\Gamma$, see \cite{5,6} for example.

We define a linear mapping $s:E\longrightarrow\bar{E}$ called \textit{a
superization} by the relation $s(e_{i,\alpha}):=\tilde{e}_{i,\pi(\alpha
)}\otimes\rho_{\alpha}$. Note that the superization mapping $s$ is not unique,
here some additional degree of freedom appears, see \cite{kk,coe} for more
details concerning with the homological description of the superization
procedure. The tensor product of spaces with $Z_{2}$ gradation is denoted by
$\otimes^{\prime}$. We have here the following relation between our tensor
products $(s\otimes^{\prime}s)(e_{i,\alpha}\otimes e_{j,\beta}):=\sigma
(\alpha,\beta)\ s(e_{i,\alpha}\otimes e_{j,\beta})$.

For the odd color Hecke $R$-matrix (\ref{ope}) there is the following
relation
\begin{equation}%
\begin{array}
[c]{l}%
(s\otimes^{\prime}s)\circ R(\alpha,\beta)=(\bar{R}(\pi(\alpha),\pi(\beta
))\hat{\otimes}R_{\delta}(\alpha,\beta))(s\otimes^{\prime}s),\label{red}%
\end{array}
\end{equation}
where $\bar{R}$ is a super-$R$-matrix on $\bar{E}$ (see \cite{13,12})
\begin{align}
&\bar{R}=\bigoplus\limits_{a,b\in Z_{2}}((-1)^{a}q^{1-2a}\Sigma_{i\in I_{a}%
}e_{i,a\ i,a}\otimes e_{i,a\ i,a}+\Sigma_{i\in I_{a},j\in I_{b},i\neq
j}(-1)^{ab}e_{j,b\ i,a}\otimes e_{i,a\ j,b}\nonumber\\&+(q+q^{-1})\Sigma_{i\in I_{a},j\in
I_{b},i<j}e_{i,a\ i,a}\otimes e_{j,b\ j,b}), \label{rsu}%
\end{align}
$R_{\delta}$ is given on the space $E^{\prime}$ by the formula (\ref{col}),
$\delta i=\pi(\alpha)$, $\delta j=\pi(\beta)$, and the tensor product $\bar
{R}\hat{\otimes}R_{\delta}$ is defined on $\tilde{E}\otimes E^{\prime}$ by the
following formula $\bar{R}\hat{\otimes}R^{\prime}=P(\bar{R}\otimes R^{\prime
})P^{-1}$, where $P:(\tilde{E}\otimes E^{\prime})^{\otimes2}\longrightarrow
\tilde{E}^{\otimes2}\otimes E^{\prime\,\otimes2}$, and $P((\tilde{e}%
_{i,\alpha}\otimes\rho_{\alpha})\otimes(\tilde{e}_{j,\beta}\otimes\rho_{\beta
}))=(\tilde{e}_{i,\alpha}\otimes\tilde{e}_{j,\beta})\otimes(\rho_{\alpha
}\otimes\rho_{\beta})$ (see \cite{8}).

Note that if we have $\Gamma_{0}=\Gamma$, then the color Hecke $R$-matrix
(\ref{red}) is even and reduces in a similar way to the ordinary $R$-matrix of
the series $A_{n-1}$, \cite{13}.

\section{DIFFERENTIAL CALCULI}

\label{diff}Here we generalize the Wess-Zumino (see \cite{17}) notion of
noncommutative differential calculi on quantum plane to the graded case. First
we generalize the concept of quantum (super-) plane to the graded case. Let
$E$ be a $\Gamma$-graded vector space equipped with a homogeneous basis
$e_{i,\alpha}$, $i\in I_{\alpha}$, $\alpha\in\Gamma$. Let $B$ be an arbitrary
graded Yang-Baxter operator on $E$ and let $\{B(\alpha,\beta):\alpha,\beta
\in\Gamma\}$ be a corresponding set of exchange operators.

An algebra defined as a quotient space ${\mathcal{A}}\equiv{\mathcal{A}%
}(E):=TE/I_{B}$, where $I_{B}$ is an ideal generated by elements
\begin{equation}%
\begin{array}
[c]{ll}%
e_{i,\alpha}\otimes e_{j,\beta}-B_{i,\alpha\ j,\beta}^{k,\beta\ l,\alpha
}\ e_{j,\beta}\otimes e_{l,\alpha}, & \text{for}\;\;i\neq j,\\
e_{i,\alpha}\otimes e_{i,\alpha},\;\;\text{if}\;\;\varepsilon_{\alpha}=-1, &
\text{for}\;\;i=j\label{idb}%
\end{array}
\end{equation}
is said to be \textit{a color (or holomorphic) quantum super-plane} for the
odd gradation and \textit{a color plane} for the even one.

The above definition means that ${\mathcal{A}}$ is an algebra generated by
$x_{i,\alpha}$, $i\in I_{\alpha}$ and relations
\begin{equation}%
\begin{array}
[c]{ll}%
x_{i,\alpha}\ x_{j,\beta}=B_{i,\alpha\ j,b}^{k,\beta\ l,\alpha}\ x_{k,\beta
}\ x_{l,\alpha}, & \text{for}\;\;i\neq j,\\
(x_{i,\alpha})^{2}=0,\;\;\text{if}\;\;\varepsilon_{\alpha}=-1, &
\text{for}\;\;i=j,\label{19}%
\end{array}
\end{equation}
where $x_{i,\alpha}:=P(e_{i,\alpha})$, $P:TE\longrightarrow{\mathcal{A}}$ is
the corresponding projection.

For the differential $d$ on a color quantum (super-) plane ${\mathcal{A}}$ we
assume the following standard properties
\begin{equation}%
\begin{array}
[c]{c}%
d=dx_{i,\alpha}\partial_{i,\alpha},\;\;d^{2}=0,\;\;d(uv)=(du)v+vdu,
\end{array}
\end{equation}
where $u,v\in E$. We assume that we have here the following exchange
relations
\begin{equation}%
\begin{array}
[c]{c}%
x_{i,\alpha}\ dx_{j,\beta}=C_{i,\alpha\ j,\beta}^{k,\beta\ l,\alpha
}dx_{k,\beta}\ x_{l,\alpha},\;\;\partial_{i,\alpha}\ \partial_{j,\beta
}=F_{i,\alpha\ j,\beta}^{k,\beta\ l,\alpha}\partial_{k,\beta}\ \partial
_{l,\alpha},\label{20}%
\end{array}
\end{equation}
and
\begin{equation}%
\begin{array}
[c]{l}%
\partial_{i,\alpha}\ x_{j,\beta}=1+(C^{-1})_{i,\alpha\ j,\beta}^{k,\beta
\ l,\alpha}x_{k,\beta}\ \partial_{l,\alpha}.\label{21}%
\end{array}
\end{equation}
for the noncommutative differentials and derivatives, respectively, $C$ and
$F$ are certain graded Yang-Baxter operators. We can calculate that there are
the following consistency conditions
\begin{equation}%
\begin{array}
[c]{l}%
\left[  (E_{12}-B_{12})\circ(E_{12}+C_{12})\right]  (\alpha,\beta)=0,\\
\left[  (E_{12}+C_{12})\circ(E_{12}-F_{12})\right]  (\alpha,\beta)=0,\\
B_{12}(\beta,\gamma)B_{23}(\alpha,\gamma)B_{12}(\alpha,\beta)=B_{23}%
(\alpha,\beta)B_{12}(\alpha,\gamma)B_{23}(\beta,\gamma),\\
C_{12}(\beta,\gamma)C_{23}(\alpha,\gamma)C_{12}(\beta,\gamma)=C_{23}%
(\alpha,\beta)C_{12}(\alpha,g)C_{23}(\beta,\gamma),\\
B_{12}(\beta,\gamma)C_{23}(\alpha,\gamma)C_{12}(\alpha,\beta)=C_{23}%
(\alpha,\beta)C_{12}(\alpha,\gamma)B_{23}(\beta,g),\\
C_{12}(\beta,\gamma)C_{23}(\alpha,\gamma)F_{12}(\alpha,\beta)=F_{23}%
(\alpha,\beta)C_{12}(\alpha,\gamma)C_{23}(\beta,\gamma).\label{con}%
\end{array}
\end{equation}
for the set of exchange operators corresponding to $B$, $C$, and $F$. The
above consistency conditions are generalization of those of Wess and Zumino
\cite{17} for the calculi on quantum plane.

The solution corresponding to the color Hecke braiding (\ref{esy}) is given by
the relations
\begin{equation}%
\begin{array}
[c]{l}%
B(\alpha,\beta)=F(\alpha,\beta)=q^{-1+2\pi(\alpha)}R(\alpha,\beta),\\
C(\alpha,\beta)=q^{1-2\pi(\alpha)}R(\alpha,\beta).\label{24}%
\end{array}
\end{equation}

We can see that this solution defines a consistent differential calculi on the
quantum space in the $\Gamma$-graded case. If we substitute the solution
(\ref{24}) into formulae (\ref{19}), (\ref{20}) and (\ref{21}), then we
obtain
\begin{equation}%
\begin{array}
[c]{l}%
x_{i,\alpha}\ x_{j,\beta}=q^{-1+2\pi(\alpha)}\varepsilon(\alpha,\beta
)\ x_{j,\beta}\ x_{i,\alpha},\;\text{for}\;\;i\neq j,\\
(x_{i,\alpha})^{2}=0,\;\;\;\text{if}\;\;\varepsilon_{i,\alpha}%
=-1,\;\;\text{and}\;\;i=j,\\
\partial_{i,\alpha}\ x_{i,\alpha}=1+q^{2(1-\pi(\alpha))}\ x_{i,\alpha
}\ \partial_{i,\alpha}+\Sigma_{k>i}(q^{2(1-\pi(\alpha))}-1)\ x_{k,\gamma
}\ \partial_{k},\\
\partial_{i,\alpha}\ x_{j,\beta}=1+q^{-1+2\pi(\alpha)}\ x_{j,\beta}%
\ \partial_{i,\alpha}.
\end{array}
\end{equation}

Next one can consider the quantum groups corresponding to the color Hecke
braiding as a new example of braided groups of Majid \cite{8}. One can also
see that the differential calculus considered here is covariant with respect
to these quantum groups.

\section{CONCLUSIONS}

Next one can generalize the deformations of the universal enveloping algebras
of the Lie algebras $A_{n-1}$ corresponding to the color Hecke braiding. The
multiparameter deformations of the universal enveloping algebras introduced in
\cite{19} can be considered as a particular example of such generalization
corresponding to the arbitrary gradation.\medskip\

\noindent\textbf{Acknowledgments}. The author would like to thank to J.
Lukierski for discussions, to A. Borowiec and R. Ra{\l }owski for critical
remarks, to K. Rapcewicz for the reading the manuscript, and to C. Juszczak
for any other help.


\begin{thebibliography}{99}                                                                                               %


\bibitem {4}{ J. Lukierski and V. Rittenberg Phys. Rev. \textbf{D18},
385 (1970). }

\bibitem {rw}{ V. Rittenberg, D. Wyler, J. Math. Phys. \textbf{19},
2193, (1978). }

\bibitem {rwu}{ V. Rittenberg, D. Wyler, N. Phys. \textbf{B139}, 189,
(1978). }

\bibitem {5}{ M. Scheunert, J. Math. Phys. \textbf{20}, 712 (1979). }

\bibitem {kn}{ Y. Kobayashi, S. Nagamachi, J. Math. Phys. \textbf{25},
3367, (1984). }

\bibitem {tro}{ R. Trostel, J. Math. Phys. \textbf{25}, 3183, (1984). }

\bibitem {kk}{ K. Kwa\'{s}niewski, Clifford and Grassmann like algebras,
Proceedings of the 13th Winter School on Abstract Analysis, Srni, 5--12
January 1987. }

\bibitem {6}{ W. Marcinek, Generalized Lie algebras and related topics,
Preprint ITP UWr No 691 and 692 (1987). }

\bibitem {WM2}{ W. Marcinek, Generalized Lie--Cartan pairs, Rep. Math.
Phys. \textbf{27}, 385 (1989). }

\bibitem {WM3}{ W. Marcinek, Graded algebras and geometry based on
Yang--Baxter operators, J. Math. Phys. \textbf{33}, 1631 (1992). }

\bibitem {10}{ W. Marcinek, On unital braidings and quantizations,
Preprint ITP UWr No 847 (1993); Rep. Math. Phys. \textbf{34}, 325, (1994). }

\bibitem {RM}{ R. Matthes, A covariant differential calculus on the
"quantum group" $C^{n}_{q}$, Proceedings of the Wigner Symposium, Goslar
(1992). }

\bibitem {9}{ V. V. Lyubashenko, Vectorsymmetries, in Seminar on
Supermanifolds, No 14, ed. by D. Leites, Stockholm University Report No 19, 1
(1987). }

\bibitem {Gur}{ D. Gurevich, Generalized translation operators on Lie
groups, Soviet. J. Contempory Math. Anal. \textbf{18} (1983). }

\bibitem {G}{ D. Gurevich, Quantum Yang--Baxter equation and a
generalization of the formal Lie theory in Seminar on Supermanifolds,
Stockholm University, Report No 19, 33 (1987). }

\bibitem {3}{ D. I. Gurevich, Soviet Math. Dokl. 38, 555-559 (1988). }

\bibitem {man}{ Y. Manin, Quantum groups and noncommutative geometry,
CRM, Montreal (1988). }

\bibitem {WM5}{ W. Marcinek, On S-Lie-Cartan pairs, in Spinors, Twistors
and Clifford algebras, ed. by Z. Oziewicz et al, Kluwer Acad. Publ. (1993). }

\bibitem {8}{ S. Majid, Braided groups, Preprint DAMTP /90-42 (1990). }

\bibitem {Maj}{ S. Majid, J. Math. Phys. \textbf{34}, 1176, (1993). }

\bibitem {bm}{ S. Majid, J. Math. Phys. \textbf{34}, 2045, (1993). }

\bibitem {qm}{ S. Majid, J. Geom. Phys. \textbf{13}, 169, (1994). }

\bibitem {HL}{ L. Hlavaty, Quantum Braided Groups, Preprint
PRA-HEP-92/18, (1992). }

\bibitem {iuz}{ K. N. Ilinski and V. M. Uzdin, Mod. Phys. Lett.
\textbf{8}, 2657, (1993). }

\bibitem {13}{ N. Yu. Reshetikin, L. A. Takhtadzhyan and L. D. Faddeev,
Leningrad Math. J. \textbf{1}, 193, (1990). }

\bibitem {12}{ D. Chang, I. Phillips and L. Rozansky, J. Math. Phys.
\textbf{33}, 3710, (1992). }

\bibitem {17}{ J. Wess and B. Zumino, Covariant differential calculus on
the quantum hyperplane, CERN-TH-5697/90, (1990). }

\bibitem {2}{ S. Majid, Int. J.Mod. Phys. A5, 1-91 (1991). }

\bibitem {7}{ Z. Oziewicz, Lie algebras for arbitrary grading group, in
Differential Geometry and its Applications, ed. by J. Janyska and D. Krupka,
World Scientific, Singapore(1990),148-152. }

\bibitem {11}{ J. Baez, Lett. Math. Phys. \textbf{23}, 133 (1991). }

\bibitem {14}{ J. C. Wallet, J. Phys. \textbf{A25}, L1159 (1992). }

\bibitem {15}{ A. Schirmacher, J. Phys. \textbf{A24}, L1249 (1991). }

\bibitem {16}{ A. Sudbery, J. Phys. \textbf{A23}, L697 (1990). }

\bibitem {18}{ S. Majid, Anyonic quantum group, Preprint DAMTP /91-16
(1991). }

\bibitem {coe}{ W. Marcinek, Colour extension of Lie algebras and
superalgebras, Preprint ITP, UWr No 746 (1990). }

\bibitem {19}{ J. F. Cornwell, J. Math. Phys. \textbf{33}, 3963 (1992).
}
\end{thebibliography}
 \end{document}